# On formula to compute primes and the $n^{th}$ prime

**Issam Kaddoura**

Lebanese International University
Faculty of Arts and Sciences, Lebanon
Email: issam.kaddoura@liu.edu.lb

**Samih Abdul-Nabi**

Lebanese International University
School of Engineering, Lebanon
Email: samih.abdulnabi@liu.edu.lb

**Abstract**

In this paper, we propose a new primality test, and then we employ this test to find a formula for $\pi$ that computes the number of primes within any interval. We finally propose a new formula that computes the nth prime number as well as the next prime for any given number.

## 1. Introduction

Since Euclid [12], primes and prime generation were a challenge of interest for number theory researchers. Primes are used in many fields, one just need to mention for example the importance of primes in networking and certificate generation [1]. Securing communication between two devices is achieved using primes since primes are the hardest to decipher [5].

The search for prime numbers is a continuous task for researchers. Some like [4] are looking for twin primes others like [11] are looking for large scaled prime numbers. The prime counting function is a function that gives the number of primes that are less than or equal to a given number.

Many like [6] and others have presented the formula to compute the number of primes between 1 and a given integer n.

This paper is divided as follow: in section 2 we present the primality test. In section 3, we introduce the prime counting function that we will use in section 5 to find the next prime to any given number. In section 4, we conduct some results and build the n<sup>th</sup> prime function. In section 6, we will use the primality test to compare our results with some recent results in the literature and conclude this paper.

## 2. Primality test

In the paper, we employ the Euclidian algorithm, Sieve of Eratosthenes and the fact that every prime is of the form 6k±1 where k an integer.

Let $x$ be a real number, the floor of $x$, denoted by $\lfloor x \rfloor$ is the largest integer that is less or equal to $x$. To test the primality of $x$ it is enough to test the divisibility of $x$ by all primes $\leq \lfloor \sqrt{x} \rfloor$.

Let $S_1(x)$ be of the form

$$S_1(x) = \frac{(-1)^{\left\lfloor \frac{\lfloor \sqrt{x} \rfloor}{6} \right\rfloor + 1}}{\left\lfloor \frac{\lfloor \sqrt{x} \rfloor}{6} \right\rfloor + 1} \sum_{k=1}^{\left\lfloor \frac{\lfloor \sqrt{x} \rfloor}{6} \right\rfloor + 1} \left\lfloor \left\lfloor \frac{x}{6k+1} \right\rfloor - \frac{x}{6k+1} \right\rfloor \qquad (2.1)$$

Similarly, let $S_2(x)$ be of the form

$$S_2(x) = \frac{(-1)^{\left\lfloor \frac{\lfloor \sqrt{x} \rfloor}{6} \right\rfloor + 1}}{\left\lfloor \frac{\lfloor \sqrt{x} \rfloor}{6} \right\rfloor + 1} \sum_{k=1}^{\left\lfloor \frac{\lfloor \sqrt{x} \rfloor}{6} \right\rfloor + 1} \left\lfloor \left\lfloor \frac{x}{6k-1} \right\rfloor - \frac{x}{6k-1} \right\rfloor \qquad (2.2)$$

**Theorem 1**: If $x$ is any integer such that g.c.d $(x, 6) = 1$
and $S(x) = \dfrac{S_1(x) + S_2(x)}{2}$, then
 (i) x is prime if and only if $S(x) = 1$
 (ii) $x$ is composite if and only if $0 \leq S(x) < 1$

**Proof:** $x$ is prime $\Rightarrow$ gcd $(x, i) = 1 \quad \forall \ 1 \leq i \leq x - 1$

$\Rightarrow x \not\equiv 0 \bmod i \ \forall \ 1 \leq i \leq x - 1 \Rightarrow \left\lfloor \left\lfloor \dfrac{x}{6k+1} \right\rfloor - \dfrac{x}{6k+1} \right\rfloor = -1$,

$\left\lfloor \left\lfloor \dfrac{x}{6k-1} \right\rfloor - \dfrac{x}{6k-1} \right\rfloor = -1 \ \forall \ k$ within the range of the summation in the formulas of $S_1(x)$ and $S_2(x)$.

$$\Rightarrow \quad S_1(x) = \frac{(-1)^{\left\lfloor \frac{\lfloor\sqrt{x}\rfloor}{6}\right\rfloor+1}}{\left\lfloor \frac{\lfloor\sqrt{x}\rfloor}{6}\right\rfloor+1} \sum_{k=1}^{\left\lfloor \frac{\lfloor\sqrt{x}\rfloor}{6}\right\rfloor+1} (-1) = 1$$

Similarly $S_2(x) = 1$ and consequently $S(x) = 1$

The proof of the second part of the theorem is obvious.

## 3. Prime Counting Function

The prime counting function, denoted by the Greek letter $\pi(n)$, is the number of primes less than or equal to a given number n. Computing the primes is one of the most fundamental problems in number theory. You can see [10] for the latest works regarding prime counting functions.

Using the previous primality test, we define the following new form of the prime counting function $\pi$.

Recall that

$$\lfloor S(i) \rfloor = \begin{cases} 1 & \text{if i is prime} \\ 0 & \text{if i is composite} \end{cases} \qquad (3.1)$$

Then $\sum_{i=m}^{n} \lfloor S(i) \rfloor$ counts the primes between m and n where $m \leq n$.

And we can write a formula for $\pi(n)$ as follows :

$$\pi(n) = 4 + \sum_{i=7}^{n} \lfloor S(i) \rfloor \qquad (3.2)$$

The size of this summation can be dramatically reduced by considering only $i$ of the form 6j+5 or 6j+7.

$$\pi(x) = 4 + \sum_{j=1}^{\left\lfloor \frac{x-1}{6} \right\rfloor} \lfloor S(6j+1) \rfloor + \sum_{j=1}^{\left\lfloor \frac{x+1}{6} \right\rfloor} \lfloor S(6j-1) \rfloor \qquad (3.3)$$

Thus the following theorem is already proved.

**Theorem 5**: $\forall x \geq 7$, $\pi(x)$ gives the number of primes $\leq x$.

## 4. The n<sup>th</sup> Prime Function

We are now ready to introduce our new formula to find the n<sup>th</sup> prime. The n<sup>th</sup> prime number is denoted by $p_n$ with $p_1 = 2, p_2 = 3, p_3 = 5$ and so on.

First we introduce $f_n(x)$ as follows

$$f_n(x) = 1 - \left\lfloor \frac{x}{n} \right\rfloor \quad \text{For } n = 1, 2, 3 \ldots \text{ and } x = 0, 1, 2 \ldots \tag{4.1}$$

Or

$$f_n(x) = \left\lfloor \frac{2n}{x+n+1} \right\rfloor \quad \text{For } n = 1, 2, 3 \ldots \text{ and } x = 0, 1, 2 \ldots \tag{4.2}$$

These functions have the property that

$$f_n(x) = \begin{cases} 1 & \text{for } x < n \\ 0 & \text{for } x \geq n \end{cases} \tag{4.3}$$

It is well known that $P_n \leq 2(\lfloor n \log n \rfloor + 1)$; see [8] and [2] for more details.

Using the following formula combined with the above formula for $\pi$

$$P_n = 7 + \sum_{x=7}^{2(\lfloor n \log n \rfloor + 1)} f_n(\pi(x)) \tag{4.4}$$

We use $f_n(x)$ as in (4.1) to obtain the following formula for $n^{th}$ prime in full:

$$P_n = 7 + \sum_{x=7}^{2(\lfloor n \log n \rfloor + 1)} \left( 1 - \left\lfloor \frac{4}{n} + \frac{1}{n} \left( \sum_{j=1}^{\lfloor \frac{x-1}{6} \rfloor} \lfloor S(6j+1) \rfloor + \sum_{j=1}^{\lfloor \frac{x+1}{6} \rfloor} \lfloor S(6j-1) \rfloor \right) \right\rfloor \right)$$

$$P_n = 3 + 2\lfloor n \log n \rfloor - \sum_{x=7}^{2(\lfloor n \log n \rfloor + 1)} \left\lfloor \frac{1}{n} \left( 4 + \sum_{j=1}^{\lfloor \frac{x-1}{6} \rfloor} \lfloor S(6j+1) \rfloor + \sum_{j=1}^{\lfloor \frac{x+1}{6} \rfloor} \lfloor S(6j-1) \rfloor \right) \right\rfloor \tag{4.5}$$

Or using $f_n(x)$ as in (4.2) to obtain the formula for $n^{th}$ prime in full:

$$P_n = 7 + \sum_{x=7}^{2(\lfloor n \log n \rfloor + 1)} \left\lfloor \frac{2n}{\pi(x) + n + 1} \right\rfloor =$$

$$7 + \sum_{x=7}^{2(\lfloor n \log n \rfloor + 1)} \left\lfloor \frac{2n}{5 + n + \sum_{j=1}^{\lfloor \frac{x-1}{6} \rfloor} \lfloor S(6j+1) \rfloor + \sum_{j=1}^{\lfloor \frac{x+1}{6} \rfloor} \lfloor S(6j-1) \rfloor} \right\rfloor \tag{4.6}$$

These formulas are in terms of n alone and we do not need to know any of the previous primes.
See [2] for formulas of the same nature.

The Wolfram Mathematica implementation of $P_n$ as in (4.5) is as follow:

```
A[x_] := (-1/(Floor[(Floor[Sqrt[x]])/6] + 1))* Sum[Floor[Floor[(x/(6 k + 1))] - (x/(6 k + 1))],
                {k, 1, Floor[(Floor[Sqrt[x]])/6] + 1}]
SB[x_] := (-1/(Floor[(Floor[Sqrt[x]])/6] + 1))* Sum[Floor[Floor[(x/(6 k - 1))] - (x/(6 k - 1))],
                {k, 1, Floor[(Floor[Sqrt[x]])/6] + 1}]
SS[x_] := (SA[x] + SB[x])/2
PN[x_] := 4 + Sum[Floor[SS[6 j + 1]], {j, 1, Floor[(x - 1)/6]}] + Sum[Floor[SS[6 j - 1]],
                {j, 1, Floor[(x + 1)/6]}]
PT[x_] := 3 + 2 (Floor[x*Log[x]]) - Sum[Floor[(1/x)*(4 +
                ( Sum[Floor[SS[6 j + 1]], {j, 1, Floor[(i - 1)/6]}] +
                Sum[Floor[SS[6 j - 1]], {j, 1, Floor[(i + 1)/6]}]))],
            {i, 7, 2 (Floor[x*Log[x]] + 1)}]
```

## 5. Next Prime

The function *nextp(n)* finds the first prime number that is greater than a given number n. As in [9] and using $S(x)$ as defined in section 2, it is clear that:

$$\prod_{x=n+1}^{n+i}(1-\lfloor S(x)\rfloor)=1 \quad \forall \ 1\leq i\leq nextp(n)-n-1$$

and

$$\prod_{x=n+1}^{n+i}(1-\lfloor S(x)\rfloor)=0 \quad \forall i \ \text{such that} \ nextp(n)-n\leq i\leq 2n$$

now consider the summation

$$\sum_{i=1}^{n}\prod_{x=n+1}^{n+i}(1-\lfloor S(x)\rfloor)=\sum_{i=1}^{nextp(n)-n}\prod_{x=n+1}^{n+i}(1-\lfloor S(x)\rfloor)+\sum_{nextp(n)-n+1}^{n}\prod_{x=n+1}^{n+i}(1-\lfloor S(x)\rfloor)$$

$$=\sum_{i=1}^{nextp(n)-n}(1)+\sum_{nextp(n)-n+1}^{n}(0)=nextp(n)-n$$

finally we obtain

$$nextp(n)=n+\sum_{i=1}^{n}\prod_{x=n+1}^{n+i}(1-\lfloor S(x)\rfloor) \qquad (5.1)$$

We used the proposed primality test to implement *nextp(n)* as follow:

1) Set $k=\left\lceil\dfrac{n-1}{6}\right\rceil$
2) Set $m=6k+1$
3) If $S(m)=1$ then go to step 8
4) Set $m=6k+5$
5) If $S(m)=1$ then go to step 8
6) $k=k+1$
7) Go to step 2
8) Output the value of *m*

The Wolfram Mathematica implementation of $nextp(n)$ is as follow:

```
A[x_] := (-1/(Floor[(Floor[Sqrt[x]])/6] + 1))* Sum[Floor[Floor[(x/(6 k + 1))] - (x/(6 k + 1))],
                    {k, 1, Floor[(Floor[Sqrt[x]])/6] + 1}]
SB[x_] := (-1/(Floor[(Floor[Sqrt[x]])/6] + 1))* Sum[Floor[Floor[(x/(6 k - 1))] - (x/(6 k - 1))],
                    {k, 1, Floor[(Floor[Sqrt[x]])/6] + 1}]
SS[x_] := (SA[x] + SB[x])/2
n=Input["Input a number:"];
k=Ceiling[(n-1)/6];
m=0;
While[True,
   m=6k+1;
   If[SS[m]==1,Break[]];
   m=6k+5;
   If[SS[m]==1,Break[]];
   k=k+1;]
```

## 6. Experimental Results

We implemented our algorithm using Wolfram Mathematica version 8. Table 1 shows the results for the n[th] prime while table 2 shows the results for the next prime. Those experimental results show the complexity of our primality test

| n[th] prime | | |
|---|---|---|
| n | P(n) | Value |
| 50 | 4.1s | 229 |
| 100 | 25.58s | 541 |
| 200 | 162.19s | 1223 |
| 250 | 286.4s | 1583 |

Table 1: nth prime

| Next prime | | |
|---|---|---|
| Next to | nextp(n) | Value |
| 10^8 | 0.04s | 100000007 |
| 10^9 | 0.187s | 1000000007 |
| 10^10 | 2.012s | 10000000019 |
| 10^11 | 1.061s | 100000000003 |
| 10^12 | 43.68s | 1000000000039 |
| 10^13 | 132.242 | 10000000000037 |

Table 2: Next prime